\documentclass{IEEEtran4PSCC}
\ifCLASSINFOpdf
\else
\fi
\usepackage[cmex10]{amsmath}
\usepackage{bbm}

\usepackage[hyphens]{url}
\usepackage{amsmath, amstext, amssymb}
\usepackage{graphicx}
\usepackage{booktabs}
\usepackage{psfrag}
\usepackage{multirow, multicol}
\usepackage{hhline}
\usepackage{eurosym}
\usepackage{color}
\usepackage{eucal,bm}

\usepackage{array}
\usepackage[latin1]{inputenc}

\begin{document}
\title{Chance Constrained Optimal Power Flow with Curtailment and Reserves from Wind Power Plants}

%% To specify the authors when (number of affiliations <= 2)
\author{
\IEEEauthorblockN{Line Roald\\ G\"oran Andersson}
\IEEEauthorblockA{Power Systems Laboratory \\
ETH Zurich\\
Zurich, Switzerland\\
\{roald, andersson\}@eeh.ee.ethz.ch}
\and
\IEEEauthorblockN{Sidhant Misra\\ Michael Chertkov\\ Scott Backhaus}
\IEEEauthorblockA{Theoretical Division and CNLS \\
Los Alamos National Laboratory\\
Los Alamos, United States\\
\{sidhant, chertkov\}@lanl.gov}
}

%% To specify the authors when (number of affiliations > 2)
% \author{\IEEEauthorblockN{Author n.1\IEEEauthorrefmark{1},
% Author n.2\IEEEauthorrefmark{2},
% Author n.3\IEEEauthorrefmark{3},
% Author n.4\IEEEauthorrefmark{3} and
% Author n.5\IEEEauthorrefmark{4}}
% \IEEEauthorblockA{\IEEEauthorrefmark{1} Department Name of Organization A\\
% Name of the organization A,
% Address A\\ Emails if wanted}
% \IEEEauthorblockA{\IEEEauthorrefmark{2} Department Name of Organization B\\
% Name of the organization B,
% Address B\\ Emails if wanted}
% \IEEEauthorblockA{\IEEEauthorrefmark{3} Department Name of Organization C\\
% Name of the organization C,
% Address C\\ Emails if wanted}
% \IEEEauthorblockA{\IEEEauthorrefmark{4}Department Name of Organization D\\
% Name of the organization D,
% Address D\\ Emails if wanted}
% }

\maketitle

\begin{abstract}
Over the past years, the share of electricity production from wind power plants has increased to significant levels in several power systems across Europe and the United States. In order to cope with the fluctuating and partially unpredictable nature of renewable energy sources, transmission system operators (TSOs) have responded by increasing their reserve capacity requirements and by requiring wind power plants to be capable of providing reserves or following active power set-point signals.
This paper addresses the issue of efficiently incorporating these new types of wind power control in the day-ahead operational planning. We review the technical requirements the wind power plants must fulfill, and propose a mathematical framework for modeling wind power control. The framework is based on an optimal power flow formulation with weighted chance constraints, which accounts for the uncertainty of wind power forecasts and allows us to limit the risk of constraint violations.
In a case study based on the IEEE 118 bus system, we use the developed method to assess the effectiveness of different types of wind power control in terms of operational cost, system security and wind power curtailment.
\end{abstract}


% Use this to place sponsorships
%\thanksto{Applicable sponsors, if any, should be placed using the \emph{thanksto} command}

\section{Introduction}
Over the last decade, electricity production from wind power plants has reached significant levels in several regions of Europe and the United States.
%Congestion management is becoming more complicated as power flow patterns change with the weather conditions and generation is shifted from traditional power plants close to the customers to more remote areas with ample wind resource, but less strong grids.
The forecast errors and fluctuations inherent to wind power generation has lead transmission system operators (TSOs) to reassess their reserve dimensioning policies \cite{CAISOreserves}. In systems with large wind penetrations, such as Denmark \cite{EnergiNet} or Ireland \cite{EirGrid}, the grid codes now require wind power plants to be able to provide active power control. These control capabilities include droop control for frequency stabilization, down-regulation of the active power output for provision of spinning reserves, capability of following an active power set-point signal and enforcing a cap on the total wind power generation \cite{EnergiNet}, \cite{Aho}.
% Maybe say something more about why these capabilities are important?

With the above mentioned capabilities, wind power plants are able to participate in ancillary service provision, system balancing and congestion management. While increased controllability is generally improving system performance, wind power plants still differ from conventional generators in that their generation output is fluctuating and their capacity is not fully known in day-ahead operational planning due to forecast uncertainty.
Therefore, if ancillary service provision from wind power plants is not planned appropriately, the use wind power control might significantly increase operational risk. For example, wind power fluctuations can render the wind power plants unable to deliver the required reserve capacities in real-time.
%Further, wind power output caps can sometimes increase wind variability, thus having a negative effect on overall system performance.

In this paper, we address the question of how to optimally incorporate the use of wind power control in day-ahead planning.
%We review the wind power control capabilities, formulate the control policies associated with the types of control we consider most relevant, and incorporate those as tractable extensions to an optimal power flow (OPF) formulation.
We account for wind power variability through the use of the weighted chance constraints (WCC) developed in \cite{CDC}. The WCC-OPF in \cite{CDC} is an alternative to the standard Chance-Constrained OPF (CC-OPF) that limits the probability of constraint violations  \cite{maria, line, 12BCH}.
Our reason for choosing the WCC-OPF is two-fold.
First, the WCCs accounts for the magnitude of constraint violation via the use of a weight function that assigns higher risk values to larger overloads. The weight function can be motivated using similar arguments as risk functions applied in risk-based OPF (e.g., \cite{Xiao2009, 14PSCC}), although the WCC limits the expected risk (as opposed to the risk of the forecasted operating point \cite{Xiao2009} or the worst-case risk with a given probability \cite{14PSCC}).
Second, while wind power curtailment as a reduction of the mean wind power production can easily the incorporated in a CC-OPF formulation, the use of WCCs enables us to model caps on the actual wind power output (e.g., to only curtail wind when the fluctuations are large) while maintaining convexity of the OPF formulation. Convexity is crucial to design efficient optimization algorithms for the WCC-OPF that scale well.

We review the wind power control capabilities, formulate the control policies associated with them, and incorporate them as tractable extensions to the WCC-OPF. In this paper, we model two types of wind control capabilities, a droop control policy that enables the wind power plants to provide reserves, and a control policy that enforces a cap on the maximum wind power output, while also accounting for the effect of the remaining uncertainty in the system. Leveraging the convexity of the WCC-OPF, we implement an efficient successive cutting-plane algorithm, and test our formulation on a modified version of the IEEE 118 bus system with 25 wind power plants. We quantify the benefits of each type of wind power control in terms of operational cost, system security and the amount of wind power curtailed.
%The WCCs ensure that the probability and extent of overloads on both generators and transmission lines remain small,

The remainder of the paper is structured as follows: Section \ref{sec:actControl} reviews the relevant control capabilities of the wind power plants. In Section \ref{sec:OPF}, the full OPF formulation is presented, with particular focus on the mathematical modeling of wind power control. The definition and handling of the weighted chance constraint is described in Section \ref{sec:WCC}.
Section \ref{sec:CaseStudy} presents the case study, and illustrates important aspects of the developed method, while Section \ref{sec:Conclusion} summarizes and concludes the paper.

\subsection{Notations}
We denote vectors by lower case letters $p, \omega$.
The components of the vectors are denoted by using subscripts, i.e, the $i$th component of $p$ is denoted by $p_i$.
Matrices are denoted by upper/lower bold case letters, $\mathbf{M}$, and $\mathbf{M}_{(i,\cdot)}, \mathbf{M}_{(\cdot,i)}$ denote the $i^{th}$ row and column of $\mathbf{M}$, respectively.
Index $i$ refers to generators, index $j$ to wind power plants, and index $ij$ refers to lines.

\section{Active power control from wind turbines}
\label{sec:actControl}
Current grid-codes in countries with significant penetration of electrical energy from wind power \cite{EnergiNet, EirGrid} require new installations of wind power plants to provide a variety of active power controls to stabilize the grid frequency and balance the system. In this paper, we consider the situation where the wind turbine adjusts the active power output by changing the amount of energy extracted from the wind (e.g., through pitch control \cite{Aho}) to follow a reference signal from the TSO. While wind power plants can provide reserves in different ways, we consider the following control mechanisms to be the most suitable for system balancing and congestion management:\\
\emph{$\Delta$P control:} The wind power plants monitors the maximum available wind power (given by the current local wind condition) and keeps the output $\Delta$MW below the maximum. The TSO can ask the wind power plant to implement this control to, e.g., curtail excess wind energy, relieve congestion, or to keep wind power capacity available for reserve provision.\\
\emph{Output cap:} The output cap is an absolute cap on the active power provided from a single wind power plant. The wind power plant produces at maximum as long as the maximum is below the output cap. If the maximum available power exceeds the output cap, the wind power output is kept constant at the cap. The TSO can use this policy to reduce the variability of the wind power output (with a low cap, the output will essentially be constant) or to handle local transmission constraints by decreasing the peak production.

The different types of control are shown schematically in Fig. 1. Notice that the power output (red line) is always less or equal to the maximum possible production (blue line) at any given point in time. In case of reserve provision, a nominal curtailment is necessary for the turbine to be able to increase the output power in reaction to a control signal from the TSO. In the case of an output cap, any production above the cap will be curtailed.
In the following, we will present a mathematical model of the two types of control, and show how they can be incorporated in an OPF formulation.

\begin{figure}[]
%%    \psfrag{116}[cc][][1][0]{116}
\includegraphics[width=0.9\columnwidth]{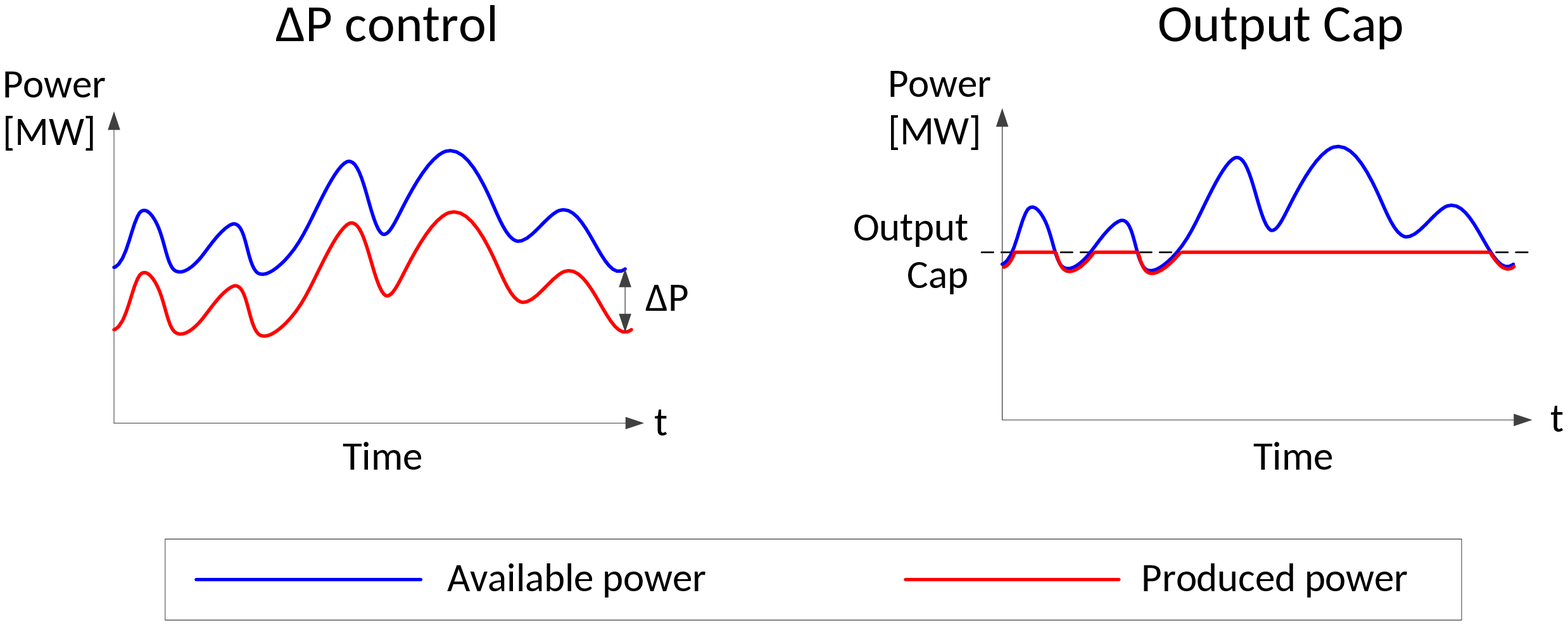}
\centering
\caption{The different types active power control ($\Delta P$ control and output caps), with wind power generation as a function of time (left) and the "controlled" wind power curve (right).}
\label{fig:actPow}
\end{figure}

\section{Optimal Power Flow with Wind Power Control}
\label{sec:OPF}
In this section, we provide a mathematical model that includes modeling of the wind power control described in the previous section within an OPF probabilistic constraints using the Weighted Chance Constraints.
The formulation extends those in \cite{CDC, 12BCH} by introducing two new modeling elements:
\begin{itemize}
	\item[(i)] the deterministic and probabilistic constraints are now modified using additional variables representing reserves from generator and wind power plants
	\item[(ii)] the generation control policies are adapted to represent the two types of wind power control capabilities.
\end{itemize}

In what follows, we will first gradually build up the objective and constraints that describe the WCC-OPF and describe each step in detail.

\subsection{Network Representation}
We represent the power transmission network as an undirected graph $\mathbb{G} = (\mathcal{V}, \mathcal{E})$, where $\mathcal{V}$ is the set of nodes with $|\mathcal{V}| = m$ and $\mathcal{E}$ is the set of edges/transmission lines of the system with $|\mathcal{E}| = n$.
The set of wind generators is denoted by $\mathcal{W}\subseteq\mathcal{V}$.
%The available wind power is considered to be a random variable, but as explained above, we assume that the wind power plants are able to control their power output to any set-point below the currently available wind power. Wind is considered as the only source of fluctuations in this paper, although the formulation can also be extended to handle fluctuations from other sources such as solar PV or load.
The set of non-wind generators is denoted by $\mathcal{G} \subseteq \mathcal{V}$, and are assumed to be controllable within their limits.
To simplify notation, we assume that each node of the system has one controllable regular generator with power output $p$, one wind generator with power output $w$ and one load with demand $d$. The nodes without generation or load can be handled by setting the respective entries to zero.

%Consider the following variables for each wind turbine $i$. The available wind power $\bar{p}_{W,i}$ is modeled as the sum of the forecasted wind power $\mu_{i}=\mathbb{E}[\bar{p}_{W,i}]$ and a zero mean fluctuating component $\omega$. The scheduled wind power production is denoted by $p_{W,i}$.

\subsection{System Modeling with Uncertain Wind Power in-feeds}

\subsubsection{Modeling Controlled Wind Power in-feeds}
The wind in-feeds $w_j$ of an uncontrolled wind power plant is modeled as the sum of the forecasted electricity production from wind, given by $\mu_j=\mathbb{E}[w_j]$, and a zero mean fluctuation around the mean, given by ${\omega_j}$, such that
\begin{align}
	w_j = \mu_j + \omega_j \label{total_wind}
\end{align}
The total wind power fluctuation is denoted by $\Omega \doteq \sum_{j\in\mathcal{W}}{\omega_j}$.
When we introduce wind power control, the output of the wind power plant changes depending on the control policy implemented.
Both the mean $\mu$ and the (distribution of the) uncertain deviation $\omega$ can be affected by the control, and we denote the new values by $v$ and $\tilde{\omega}$ respectively. The new total wind deviation from mean is given by $\tilde{\Omega} = \sum_{j} \tilde{\omega}_j$.

\emph{Wind Power Plants with $\Delta P$ control}
The TSO might use $\Delta P$ control to reduce the scheduled mean power from a wind power plant $j \in \mathcal{W}$ to a value $v_j \leq \mu$.
%In this case, the random fluctuation component $\omega_i$ is unaffected.
The wind power output of such plants are given by
 \begin{align}
 	\tilde{w}_j = v_j + \omega_j. \label{wind_v}
 \end{align}
%The wind power plants with $\Delta P$ control can use the unused part of the energy to provide reserves. We assume that these wind power plants provide reserves by following an affine control policy which changes their output by $- \alpha_i \tilde{\Omega}$. The resulting wind power output of this plant is given by
We assume that the wind power fluctuations $\tilde{\Omega}$ are balanced through the automatic generation control (AGC), with participation factors $\alpha$. For wind power plants contributing to AGC, the resulting wind power output is thus given by
 \begin{align}
 	\tilde{w}_j = v_j - \alpha_j \tilde{\Omega} + \omega_j. 	\label{deltap_output}
 \end{align}
We denote the set of wind power plants providing reserves by $\mathcal{R}\subset\mathcal{C}$.
%In order to reduce congestion or accommodate reserve provision, the TSO may reduce the scheduled wind power output to a value below $\mu$. We denote this mean scheduled wind power output by $v$. Note that when the TSO reduces the scheduled output to $v<\mu$, the wind power fluctuations $\omega$ are not affected (they just fluctuate around a new mean value).
%Further, the TSO is allowed to introduce output caps on the wind power plants, denoted by $\bar{\omega}$. The output cap is defined relative to the scheduled output $v$, such that a negative cap implies a cap below $v$ and a positive cap implies a cap above $v$. In this case, the distribution of the wind power fluctuation is truncated above $\bar{\omega}$.
%Note that $\bar{\omega}$ are not subject to optimization, but are given as a pre-defined input.

\emph{Wind power plants with Output Caps}
With output cap control, the TSO can enforce an upper limit on the power output of individual wind power plants. In this case, the power output is given by
\begin{align}
	\tilde{w}_j = v_j + \min\{ \omega_j, \bar{\omega}_j \}, 		\label{cap_output}
\end{align}
where $v_j \leq \mu_j$ is the mean power output and the cap $\bar{\omega}$ is the upper bound on the wind power fluctuation. We denote the set of wind power plants with output caps by $\mathcal{C}\subset\mathcal{W}$, and assume that these wind power plants do not provide reserves, i.e., $\alpha_j = 0~\forall_{j \in \mathcal{C}}$.
%When one or more wind power plants are subject to output caps, the total wind power fluctuation $\Omega$ needs to be redefined. When any of the wind fluctuations exceed their output cap, the output will be constant, and we define the new $\tilde{\Omega}$ as
%\begin{equation}
%\tilde{\Omega} = \sum_{j\in\mathcal{W}/\mathcal{C}} \omega_j + \sum_{k\in\mathcal{C}} \min(\omega_k, \bar{\omega}_k)~,
%\end{equation}
%where $\mathcal{C}$ is the set of wind power plants with output caps.
%For simplicity, we assume that the output caps $\bar{\omega}_j$ are defined separately for each wind power plant, and that wind power plants with output caps are not providing reserves (i.e., $\alpha_j=0$ for wind power plants with a cap). Assuming that a wind power plant with no output cap is equivalent to a wind turbine with output cap $\bar{\omega}=\infty$, we redefine the wind power outputs as
%\begin{equation}
%\tilde{v}_j = v_j + \min(\omega_j,\bar{\omega}_j) - \alpha_j \tilde{\Omega}~, \quad \forall_{j\in\mathcal{W}}
%\end{equation}

\subsubsection{Modeling Generation Output}
We assume that the regular generators participate in reserve provision, and their power output is given by
\begin{align}
	\tilde{p}_i = p_i - \alpha_i \tilde{\Omega},  \label{gen_output}
\end{align}
where $p_i$ is the nominal power output set point and $\alpha_i$ is the participation factor in balancing the wind power fluctuations.

\subsubsection{Modeling Power Flows}
The power flows on each transmission line is computed according to the standard DC approximation \cite{wollenberg},
\begin{equation}
p_{ij}=\mathbf{M}_{(ij,\cdot)}( \tilde{p} + \tilde{w}-d). \forall_{ij\in \mathcal{E}}.
\label{eq:lineflows}
\end{equation}
The matrix $\mathbf{M} \in \mathbb{R}^{n \times m}$ relates the line flows to the nodal power injections, which are expressed as the sum of controlled generation output $\tilde{p}$, controlled wind power production $\tilde{w}$ and demand $-d$.
The matrix $\mathbf{M}$ is defined as
\begin{equation}
\mathbf{M} = B_f \begin{bmatrix} (\widetilde{B}_{bus})^{-1}   ~~~ \bold{0}  \\ ~~\bold{0}  ~~~~~~~~~ 0\end{bmatrix}
\end{equation}
where ${B}_{f}$ is the line susceptance matrix and $\widetilde{B}_{bus}$ the bus susceptance matrix (without the last column and row) \cite{maria}. We have used
$\mathbf{M}_{(ij,\cdot)}$ as the row of $\mathbf{M}$ related to the line $(ij) \in \mathcal{E}$.

\subsection{Objective function}
The objective of the WCC-OPF is to minimize the sum of the generation and reserve cost. This is expressed as
\begin{equation}
%\min\limits_{p,v,\alpha,r^+,r^-} \ &\sum_{i\in\mathcal{G}} \left(c_i p_i + c_i^+ r_i^+ + c_i^- r_i^-\right)  \label{gen_cost} \\
%                                + &\sum_{j\in\mathcal{W}} \left(c_j p_j + c_j^+ r_j^+ + c_j^- r_j^-\right) \label{CC-OPF}
\min\limits_{p,v,\alpha,r^+,r^-} \ \sum_{i\in\mathcal{G}} c_i p_i + \sum_{j\in\mathcal{W}} c_j v_j
                                    + \sum_{i\in\mathcal{W,G}} \left(c_i^+ r_i^+ + c_i^- r_i^-\right)
\end{equation}
The vectors $c,~c^+,~c^-$ denote the cost, i.e. bids, from the generators and wind power plants for energy, up- and down reserves. The up- and down-reserves $r^+,~r^-$ are defined as non-negative,
\begin{equation}
0 \leq r^+ \leq r^{+}_{max}~, \quad 0 \leq r^- \leq r^-_{max}~.
\label{CC-OPF}
\end{equation}
and $r^+_{max},~r^-_{max}$ are upper limits that define the ability or willingness of the generators and wind power plants to provide reserves.

%The corresponding power flow constraints are given by
%\begin{alignat}{1}
%& \quad \mbox{Prob}_{\bm\omega} \left( \mathbf{M}_{(ij,\cdot)}(p-\alpha\tilde{\Omega}+\tilde{v}-d) > p_{ij}^{max} \right) < \epsilon_{ij}~, \quad  \forall_{ij \in \mathcal{E}}~, \label{chance_line_max} \\
%& \quad \mbox{Prob}_{\bm\omega} \left( \mathbf{M}_{(ij,\cdot)}(p-\alpha\tilde{\Omega}+\tilde{v}-d) < -p_{ij}^{max} \right) < \epsilon_{ij}~, \quad  \forall_{ij \in \mathcal{E}}~, \label{chance_line_min}
%\end{alignat}
%which are also probabilistic constraints with risk limit $\epsilon_{ij}$.

\subsection{Power Balance Constraints}
%The Power Flow constraints in Eq.~(\ref{eq:lineflows} )must be appended by the total power balance constraint which says that the total injection and withdrawal into the system must sum to zero,
To ensure power balance, we enforce the constraint
\begin{align}
	\sum_i \tilde{p}_i + \tilde{w}_i - d_i = 0.  \label{random_power_balance}
\end{align}
Since quantities $\tilde{p}_i$ and $\tilde{w}_i$ are random quantities whose values depend on the random wind fluctuation $\omega$, the above relation must hold for all possible realizations of $\omega$.
This can be enforced by separating the nominal part (when the fluctuation is zero) and the stochastic part of \eqref{random_power_balance}.
%The nominal power balance is obtained by setting the wind fluctuation to zero.
The nominal constraint is obtained by substituting $\omega = 0$ in Eqs.~(\ref{deltap_output}),(\ref{cap_output}),(\ref{gen_output}) and plugging in Eq.~(\ref{random_power_balance}),
\begin{align}
	\sum_{i \in \cal{G}} p_i - \sum_{i \in D} d_i + \sum_{j \in \cal{W}} v_j + \sum_{j \in \cal{C}}\min\{0,\bar{\omega}_j\} = 0. \label{nomial_powerbalance}
\end{align}
When the wind fluctuation $\omega$ is non-zero, we must have
\begin{align}
	0 &= \sum_{i \in \cal{G}} \left(p_i - \alpha_i \tilde{\Omega}\right) - \sum_{i \in \cal{D}} d_i + \sum_{j \in \cal{W}}\left(v_j +\omega_j\right)\nonumber\\
      &  + \sum_{j \in \cal{R}}\alpha_j \tilde{\Omega}+ \sum_{j \in \cal{C}}\min\{ \omega_j, \bar{\omega}_j \} \nonumber \\
%	& = \sum_{i \in \cal{G}} p_i - \sum_{i \in D} d_i + \sum_{j \in \cal{W}} v_j + \sum_{j\in\cal{C}}\min\{0,\bar{\omega}_j\}  \nonumber\\
%	 &=  \sum_{i \in \cal{G},\cal{R}} - \alpha_i \tilde{\Omega} + \sum_{j \in \cal{R}} \omega_j + \sum_{j\in\cal{C}}\min\{ \omega_j,\bar{\omega}_j\}\nonumber\\
	 &= (1 - \sum_{i} \alpha_i ) \tilde{\Omega}.
\end{align}
From the last line it follows that to ensure power balance for every value of random wind fluctuation, it is enough to enforce the constraint
\begin{align}	
	\sum_{i} \alpha_i = 1.
\end{align}

\subsection{Constraints for Line Flow Limits, Generation Limits and Reserves}
What remains is to enforce generation and transmission constraints. For quantities that are functions of fluctuating wind in-feeds, we use WCCs to ensure that \emph{(a)} the risk of having procured too little generation capacity, \emph{(b)} the risk of wind power plants not being able to provide the contracted reserves and \emph{(c)} the risk of line power flows exceeding the transmission limits, are all small in a probabilistic sense.
%In each of the above cases it is convenient to separate the deterministic and probabilistic part of these constraints. We use weighted chance constraints to enforce the probabilistic constraints.
For the sake of readability, in this section we will suppress the details of the WCCs by using a short-hand notation, and describe them in detail in the next section. More specifically, whenever we need to enforce that the risk of the quantity $y(\tilde{\omega})$ (representing e.g. a line overload ) is small, we will denote it by
\begin{align}
 WCC\left(y(\tilde{\omega}) \leq 0 \right) \leq \epsilon~,
	%y \stackrel{\tiny{wcc}}{\leq}0.
\label{eq:WCC_generalform}
\end{align}
where $\epsilon$ represents the risk limit.

\subsubsection{Constraints for Conventional Generators}
For a conventional generator, we enforce generation limits and constraints on reserve availability in the following way:
\begin{alignat}{1}
 \quad & p + r^{+}  \leq  {p}_G^{max}~, \label{gen_max2} \\
 \quad &p - r^{-}  \geq {p}_G^{min}~, \label{gen_min2} \\
 \quad &WCC\left(- \alpha_i \tilde{\Omega} {>} r^+_i \right) < \epsilon_i ~, \quad  \forall_{i \in \mathcal{G}}~, \label{chance_gen_min2} \\
 \quad &WCC\left(- \alpha_i \tilde{\Omega} {<} r^-_i \right) < \epsilon_i ~, \quad  \forall_{i \in \mathcal{G}}~.\label{chance_gen_max2}
% \quad  - &\alpha_i \tilde{\Omega} \stackrel{\tiny{wcc}}{>} r^+_i ~, \quad  \forall_{i \in \mathcal{G}}~, \label{chance_gen_min2} \\
% \quad - &\alpha_i \tilde{\Omega}  \stackrel{\tiny{wcc}}{<}  r^-_i   ~, \quad  \forall_{i \in \mathcal{G}}~.\label{chance_gen_max2}
\end{alignat}
Here, \eqref{gen_max2}, \eqref{gen_min2} %enforces the generation limits.%, with $p$ denoting the scheduled generator output, and $r^+,~r^-$ the up- and down-reserves.
enforces that the generators do not commit to providing energy and reserves that will cause them to exceed their physical minimum and maximum generation limits.
Constraints \eqref{chance_gen_min2}, \eqref{chance_gen_max2} describe the activation of reserves in reaction to the wind fluctuation $\tilde{\Omega}$.%, and %$\alpha_i$ is the fraction of the wind fluctuation that is balanced by the $i$th generator. They
They use WCCs to enforce that the risk of not having enough reserves available to cover the fluctuations remains below the risk limit $\epsilon_i$.

\subsubsection{Constraints for Wind Power Plants with $\Delta P$ Control}
The wind power plants use $\Delta P$ control to maintain constant capacities $r^+,~r^-$ available for reserves, similar to a conventional generator.
However, for the wind power plants, we also need to ensure that the risk of the wind power plant not being able to provide the promised reserves is small, since the power output of a wind power plant is dependent on the wind realization.
%However, for a wind power plant, we not only need to ensure that the reserves will be enough to cover the wind power imbalances, but also need to ensure that the risk of the wind power plant not being able to provide the promised reserves is small. The latter constraint is applicable only to wind power plants because unlike conventional generators the maximum power output capacity of a wind power plant is dependent on the wind realization.
This can be achieved by enforcing the following chance constraints:
\begin{alignat}{1}
 \quad &v + r^{+} \leq \mu~, \label{wind_max2} \\
 \quad &WCC\left(v_j + \omega_j - r_j^{-} \leq 0 \right) < \epsilon_r~, \quad  \forall_{j \in \mathcal{W}}~, \label{wind_min2} \\
 \quad &WCC\left(- \alpha_j \tilde{\Omega} {>} r^+_j \right) < \epsilon_j ~, \quad  \forall_{j \in \mathcal{W}}~, \label{chance_wind_min} \\
 \quad &WCC\left(- \alpha_j \tilde{\Omega} {<} r^-_j \right) < \epsilon_j ~, \quad  \forall_{j \in \mathcal{W}}~.\label{chance_wind_max}
% \quad - &\alpha_j \tilde{\Omega} \stackrel{\tiny{wcc}}{>} r^+_j ~, \quad  \forall_{j \in \mathcal{W}}~, \label{chance_wind_min} \\
% \quad - &\alpha_j \tilde{\Omega} \stackrel{\tiny{wcc}}{<} r^-_j ~, \quad  \forall_{j \in \mathcal{W}}~.\label{chance_wind_max}
\end{alignat}
Eq. \eqref{wind_max2} is the nominal constraint that ensures that the scheduled power generation $v$ and the capacity allocated for up-reserves $r^+$ remain below the forecasted power $\mu$. This is a deterministic constraint, since both the actual produced power and the available power will vary by $\omega$, and $\omega$ therefore cancels out.
Eq. \eqref{wind_min2} ensures that the actual produced power $v_j+\omega_j$ is high enough to provide the expected down-reserves $r^+$. As above, this equation is formulated as a WCC, which limits the risk $\epsilon_r$ that the wind power plant will not be able to provide the allocated reserve capacity.
Constraints \eqref{chance_wind_min}, \eqref{chance_wind_max} have the same interpretation as the corresponding constraints for the conventional generators \eqref{chance_gen_min2}, \eqref{chance_gen_max2}.

Note that the wind power plants are contracted to provide a constant amount of up and down reserves $r^+,~r^-$, i.e., the amount of provided reserves is not allowed to change with the available wind power. To provide up-reserves $r^+$, the wind power plants thus curtails a constant amount 1 MW wind power to provide 1 MW up-reserves.
%As noted in Section \ref{sec:actControl}, using wind power to provide reserves induces a certain amount of curtailment. With $\Delta P$ control, the turbine curtails a constant amount of $\mu - v$ MW to provide up reserves (i.e., it is necessary to curtail 1 MW of wind power to obtain 1 MW of reserves). The total curtailed power is thus $\mu - v$.
%
%The requirements/confidence for reserve availability can be adjusted by choosing an appropriate $\epsilon_{R,i}$. Choosing a higher $\epsilon_{R,i}$ impacts how much reserves each wind turbine is able to provide, but does not influence the amount of curtailment per MW of reserves.

\subsubsection{Constraints for Line Flow Limits}
To limit the risk of transmission line overloads, we enforce the power flow constraints using WCCs, i.e.,
%We enforce probabilistic WCC based constraints that the line power flows remain below their maximum.
%The corresponding power flow constraints are given by
\begin{alignat}{1}
& WCC\left(\mathbf{M}_{(ij,\cdot)}(p-\alpha\tilde{\Omega}+\tilde{v}-d) {>} p_{ij}^{max}  \right) < \epsilon_{ij}~, \quad  \forall_{ij \in \mathcal{E}}~, \label{chance_line_max} \\
& WCC\left(\mathbf{M}_{(ij,\cdot)}(p-\alpha\tilde{\Omega}+\tilde{v}-d) {<} -p_{ij}^{max} \right) < \epsilon_{ij}~, \quad  \forall_{ij \in \mathcal{E}}~, \label{chance_line_min}
%& \mathbf{M}_{(ij,\cdot)}(p-\alpha\tilde{\Omega}+\tilde{v}-d) \stackrel{\tiny{wcc}}{>} p_{ij}^{max} ~, \quad  \forall_{ij \in \mathcal{E}}~, \label{chance_line_max} \\
%& \mathbf{M}_{(ij,\cdot)}(p-\alpha\tilde{\Omega}+\tilde{v}-d) \stackrel{\tiny{wcc}}{<} -p_{ij}^{max} ~, \quad  \forall_{ij \in \mathcal{E}}~, \label{chance_line_min}
\end{alignat}
where $\epsilon_{ij}$ is the risk limit for line overloads and $p_{ij}^{max}$ is the maximum transmission capacity of line $ij$.

%\subsection{Generation and Reserve Constraints}
%The nominal power balance without any wind power fluctuations is given as
%\begin{equation}
%\sum_{i\in\cal{V}}p_i - d_i + v_i + \min(0,\bar{\omega}_i) = 0~, \label{powerbalance}
%\end{equation}
%where the last term accounts for the effect of wind power caps below the mean output.
%Since secure operation of the power system requires balance between produced and consumed power at all times, any short-term deviation $\omega$ in the wind power production must be balanced by activation of reserves. As in the previous work \cite{12BCH}, [others?], these adjustments are modeled through an affine policy, reflecting the automatic generation control which is establishing balance within tens of second to a few minutes \cite{wollenberg}.

\section{Weighted Chance Constraints}
\label{sec:WCC}
%\textbf{I would like to streamline the first part of this section, and try to focus it more on the linear WCC (e.g., explain the interpretation as an expected overload), but I can do that tomorrow.}
In this section we describe the specifics of the WCCs that should replace the constraints of the form \eqref{eq:WCC_generalform} in the previous section.
%As mentioned in the introduction, the weighted chance constraints introduced in \cite{CDC} have some advantages over the standard chance constraints. First, they allow for a more flexible definition of risk which can account for the magnitude of constraint violation. By using a weight function it is possible to assign higher risk to constraint violations of larger magnitude, e.g., by adopting weight functions similar to the risk functions in risk-based OPF \cite{Xiao2006, 14PSCC}.
%Second, the resulting WCC-OPF is a convex optimization problem even for non-affine control policies, and allow for efficient optimization algorithms.
We will first give a general introduction to the WCCs and describe the physical interpretation of a WCC with a linear weight function, which is adopted here. We then present analytical expressions for the WCCs with $\Delta P$ and output cap control under assumption of a normal distribution.
%The weighted chance constraint differs from the standard chance constraints in some important aspects.

\subsection{General Weighted Chance Constraint}
The general weighted chance constraint is a constraint of the form
\begin{equation}
\int_{-\infty}^{\infty}f(y(\omega))P(\omega)d\omega \leq \epsilon~,
\label{cvar_general_1}
\end{equation}
where $P(\omega)$ is the multivariate distribution function of the fluctuations. The quantity $y(\omega)$ denotes the magnitude of constraint violation, and is defined differently for each type of constraint, and $f(.)$ is the risk weighting function.
For a violation of the upper or lower reserve limits \eqref{chance_gen_min2}, \eqref{chance_gen_max2},  \eqref{chance_wind_min} and \eqref{chance_wind_max} the magnitude of constraint violation is given by
\begin{equation}
y(\omega) = -\alpha_i \tilde{\Omega} - r^+~, \quad y(\omega) = r^- + \alpha_i \tilde{\Omega}~, \forall_{i \in \mathcal{G, W}}~. \label{y_gen}
\end{equation}
Similarly, violations of the availability of down reserves \eqref{wind_min2} are defined by
\begin{equation}
y(\omega) = - v_j - \omega_j + r_j^{-} , \forall_{i \in \mathcal{W}}~,
\end{equation}
and violations of the upper and lower line flow limits by
\begin{equation}
y(\omega) = \tilde{p}_{ij}(\omega) - p_{ij}^{max}~, \quad 	y(\omega) = p_{ij}^{min} - \tilde{p}_{ij}(\omega)~, \forall_{ij \in \mathcal{E}}~. \label{y_line_l}
\end{equation}

Whenever we have $y>0$, it indicates a violation of the limit, while $y<0$ implies that we are in a safe operating region. The weighting function $f(y(\omega))$, which is nonzero only if $y>0$, describes the risk related to the overload, and can be chosen in different ways.
For example, \eqref{cvar_general_1} is equivalent to a standard chance constraint if $f(y)$ is the unit step function, i.e., 0 for $y<0$ and 1 for $y\geq 0$. However, the step function is not convex, which means that the standard chance constraint will not always be a convex constraint.
On the other hand, as proven in \cite{CDC}, the constraint (\ref{cvar_general_1}) is a convex for general generation control policies and general probability distributions of the wind fluctuations whenever the weight function $f(.)$ is convex.

Using a convex risk function, which assigns a higher risk to constraint violations of larger magnitude, also makes sense from an engineering point of view, and has been applied in risk-based OPF (e.g. \cite{Xiao2009}, \cite{Gabi2012}, \cite{14PSCC}). Note that the choice of the risk function affects the interpretation of the risk limit $\epsilon$. For a weighted chance constraint with a linear weight function, which is applied here, the risk limit can can be interpreted as an upper bound on the expected constraint violation. The unit of $\epsilon$ in Eq.~(\ref{cvar_general_1}) is the same as for the LHS of the constraint, i.e., we define $\epsilon$ as an acceptable expected magnitude of overload in MW.

One main advantage of  using a convex weighting function $f(y)$ in (\ref{cvar_general_1})  is the ability to handle more general control policies while still maintaining convexity, as shown in \cite{CDC}.
This is especially important when modeling wind power control with generation cap, since the associated policy is inherently non-affine and breaks the convexity of the standard chance constraints.

\subsection{Expressions for the Weighted Chance Constraints with Linear Weight Functions with and without Cap Control}
The expressions we derive for the WCC with a linear weight function assume that the line flows and generation outputs (which are weighted sums of the wind fluctuations omega) are normally distributed. While this is a relatively strong assumption, it might be justified using the Central Limit Theorem \cite{dasgupta} in systems with a large number of wind power plants. However, the WCC can also be computed even when the distribution is not normal, although this might require more complex numerical methods for evaluation of the integrals.
\\
\subsubsection{Without Cap Control}
When the wind fluctuation $\omega$ is a multivariate Gaussian random variable and there are no wind power plants enforcing cap control, then for each of the probabilistic constraints, the overload random variable  $y(\omega)$ is normally distributed. In this case, the linear weighted chance constraint can be computed by using the expression for the expectation of a truncated normal random variable:
 \begin{align}
 	& \int_{0}^{\infty} y P(y) dy = \nonumber \\
    & \mu_y \left(1 - \Phi\left( \frac{-\mu_y}{\sigma_y} \right)\right) + \frac{\sigma_y}{\sqrt{2\pi}}e^{-\frac{1}{2}\left(\frac{-\mu_y}{\sigma_y}\right)^2} \leq \epsilon~, \label{lin}
 \end{align}
 and $\mu_y$ and $\sigma_y$ denote the mean and standard deviation of $y$.
\\
\subsubsection{With Cap Control}
% When some of the wind power plants are enforcing cap control, the  resulting policy is non-affine and the computation of the linear weighted chance constraint becomes more involved. Let
When some of the wind power plants are enforcing cap control, the distribution of $y$ is no longer normal, since the standard deviation changes when an output cap is reached. The computation of the linear weighted chance constraint thus becomes more involved.
Let $\mathcal{C} = \{ i_1, \ldots, i_K \}$ denote the indices of the wind power plants implementing cap control.
For each of these wind power plants, we split the integral over ${\omega}_{\mathcal{C}}$ in \eqref{cvar_general_1} into two parts, one where the fluctuation in the wind power output is below the cap  $\bar{\omega}_{\mathcal{C}}$, and one where the fluctuation is above the cap, by defining the sets $\mathcal{S}_{k,b} \in \mathbbm{R}$ for $k = 1, \ldots, K$ and $b \in \{ 0,1\}$ as
 \begin{align*}
 	\mathcal{S}_{k,b} = \begin{cases} (-\infty, \bar{\omega}_k], \ &\mbox{for } b = 0, \\
								( \bar{\omega}_k, \infty) \ &\mbox{for } b =1.
								\end{cases}
 \end{align*}

Then, the linear weighted chance constraints can be computed as
 \begin{align}
 	&\int_{0}^{\infty}\! \int_{\!-\!\infty}^{\infty} \!  y P(y, \omega_{\mathcal{C}}) d \omega_{\mathcal{C}} dy  \label{wcc_plain} \\
    =&\sum_{b \in \{ 0,1\}^K}  \int_{0}^{\infty} \! \int_{\mathcal{S}_{1,b_1}} \ldots \int_{\mathcal{S}_{K,b_K}} \! y P (y, \omega_{\mathcal{C}}) d \omega_{\mathcal{C}} d y \leq \epsilon,  \label{wcc_specific}
 \end{align}
where $\omega_{\mathcal{C}} = ( \omega_{i_1}, \ldots, \omega_{i_K})$.  In the above summation, whenever  $b_k = 1$, the corresponding $\omega_{i_k}$ is larger than
 $\bar{\omega}_{i_k}$ and the power output is fixed at  $\bar{\omega}_{i_k}$. The overload $y$ is jointly distributed along with $\omega_\mathcal{C}$ according to a multivariate
 Gaussian distribution. One can evaluate the integrals using  efficient numerical schemes for Gaussian integration \cite{Genz95}. For smaller values of $K$, it is efficient to split (\ref{wcc_plain}) into the summation in (\ref{wcc_specific}) and evaluate each of the terms using Gaussian quadrature integration specialized to rectangular domains. When $K$ is large, it is more appropriate to perform numerical integration  directly on Eq.~(\ref{wcc_plain}) using Monte-Carlo sampling.

 We remark here that the constraints in Eq.~(\ref{wcc_specific}) are convex w.r.t. the optimization variables in the WCC-OPF, namely $p,\alpha$ and $r$, whereas the same can be shown to be non-convex for the case of standard chance constraints.

\section{Numerical Results}
\label{sec:CaseStudy}

\subsection{Implementation of the WCC-OPF with Wind Power Control}
We leverage the convexity of the WCCs to devise an efficient outer-approximation-cutting-plane algorithm, similar to the one in \cite{12BCH}. The algorithm is implemented it in Julia language \cite{julia} using JuMP \cite{jump}. At each step, a Linear Program (LP) corresponding to an outer approximation of the feasible set is solved.
The outer approximation is progressively tightened by adding tangential cutting-planes to the convex WCCs that eliminate the current infeasible solution. Each successive Linear Program is warm started using the solution of the previous one.

 Since the constraints and gradients for the WCC-OPF with cap control must be evaluated using numerical integration methods, it is desirable to reduce the number of function evaluations corresponding to these integrals. The cutting-planes algorithm accomplishes this by selectively evaluating the gradients only for the violated constraints, and shifting most of the computations over to the Linear Programming Solver.

\subsection{Case Study}

We base our study on the IEEE 118-bus system as provided with Matpower 4.1 \cite{matpower}, with a few modifications as follows. For the generator bids for energy and reserves $c,~c^+,~c^-$, we use the linear cost coefficients. For wind power plants, we assume zero marginal cost, and set $c_j=c^+_j=c^-_j=0$. Although the formulation could be extended to include unit commitment, it is not considered here.  Therefore, the minimum generation output of the conventional generators is set to zero. To obtain a more stressed system state, we increase the load by a factor of 1.25 and descrease the transmission limits by 0.75.
Wind power plants are located at 25 different buses throughout the system.
 %\textbf{we could add a picture of the system and the covariance matrix if space is available}.
Their locations, forecasted power output and correlation matrix can be found in \cite{118BusWindData}. The standard deviation of each wind power plant is set to 10\% of the forecasted power output.
When considering different levels of wind power penetration, both the forecasted power output and the standard deviations are scaled by a factor of
\begin{equation}
\frac{\sum_{i\in\mathcal{D}}{d_i}}{\sum_{j\in\mathcal{W}}{\mu_j}}\cdot \frac{X}{100}~,
\end{equation}
where $X=\{25,~50,~75,~100,~125\}$ denotes the percentage of forecasted wind power relative to the total system load.
The risk limits for all WCCs were set to $\epsilon=0.1$ MW. With the cutting-plane algorithm described above, a solution to the WCC-OPF for the 118 bus system is obtained within a couple of minutes on a laptop.

%Description of the test system (IEEE 118 bus)\\
%\textbf{Need to define what is expected wind energy curtailment}

\subsubsection{Impact of Reserve Provision from Wind Power Plants}
\begin{figure}
\includegraphics[width=0.9\columnwidth]{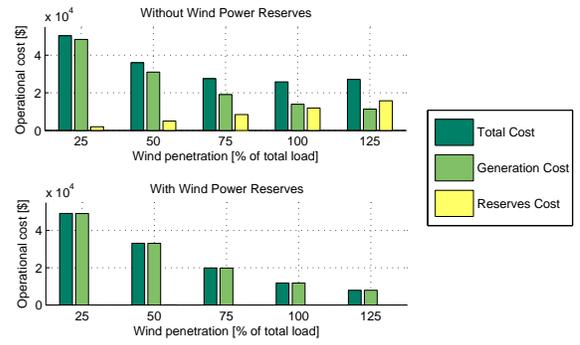}
\centering
\caption{Total cost, cost of production and cost of reserves for the two cases with and without reserves from wind power.}
\label{fig_comparison_cost}
\end{figure}

\begin{figure}
\includegraphics[width=0.9\columnwidth]{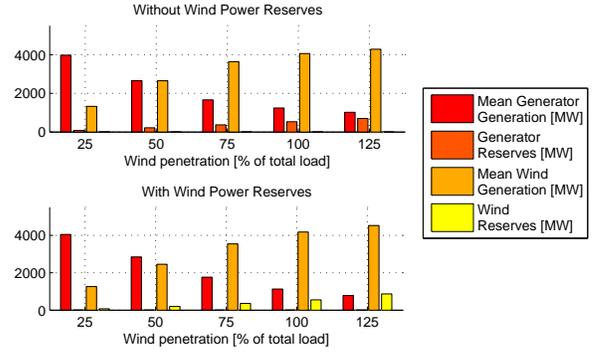}
\centering
\caption{Production and reserves from generators and wind power plants for each level of wind penetration.}
\label{fig_comparison_wind}
\end{figure}

We compare the cases when the wind power plants using $\Delta P$ control can (a) only curtail their mean as in Eq. \eqref{wind_v} and (b) the curtailed energy can be used to provide reserves as in Eq. \eqref{deltap_output}. Figure~\ref{fig_comparison_cost} shows the total cost, cost of generation and cost of reserves for the two cases without and with wind power reserves under various levels of wind penetration, and Figure~\ref{fig_comparison_wind} shows the total amount of generation and reserves provided by both conventional generators and wind power plants.

As one might expect, the total cost is lower when the wind power plants are providing reserves. This difference is especially amplified when the system
has higher wind penetration, and a large fraction of the cost difference is attributed to a higher cost of reserves. With high wind penetration, there is significantly more uncertainty in the system, and without (zero cost) wind power reserves, a large amount of reserves must be bought from conventional generators to balance fluctuations.

Additionally, in the case without wind power reserves, the nominal generation cost is also higher. 
%The generators must run at a higher nominal set-point $p$ in order to be able to provide enough down reserve $r^{-}$. 
As observed in Figure~\ref{fig_comparison_wind}, the generators must run at a higher nominal set-point $p$ in order to be able to provide enough down reserve $r^{-}$, particularly during high wind penetrations. Thus, the wind power nominal set points $v$ must be lowered, leading to under-utilization of available cheap wind power.
%In Figure~\ref{fig_comparison_wind}, this becomes apparent by noticing that the expected wind power curtailed (yellow bar) is much higher when there are no wind reserves, especially for the $100$ and $125$ percent  wind penetration cases.

%Show the value of wind power reserves (with increasing penetration)

%\begin{table}
%\caption{Generation and Reserves from Wind Power Plants and Conventional Geneators for Different Wind Penetrations $W\%$.}
%\footnotesize
%\centering
%\resizebox{0.5\textwidth}{!}{%
%\begin{tabular}{crrrrr}
%\toprule
% & \multicolumn{4}{c}{Generators} & \multicolumn{4}{c}{Wind Power Plants}  \\
% \cmidrule(lr){2-5}
% \cmidrule(lr){6-9}
% \cmidrule(lr){8-10}
% & \multicolumn{2}{c}{ $\sum{p}$ [MW] } & \multicolumn{2}{c}{\sum{r^+ + r^-}$ [MW] }  & \multicolumn{2}{c}{ $\sum{p}$ [MW] } & \multicolumn{2}{c}{\sum{r^+ + r^-}$ [MW] }\\
% $W\%$ & Without & With & Without & With & Without & With & Without & With \\
% \midrule
%  25   &3.9769 & 4.0470 &0.0816 & 0.0016 & 1.3256 & 1.2555 &-&   0.0701\\
%  50   &2.6512 & 2.8489 &0.2213 & 0.0036 & 2.6512 & 2.4536 &-&   0.1977\\
%  75   &1.6622 & 1.7617 &0.3687 & 0.0023 & 3.6403 & 3.5408 &-&   0.3557\\
% 100   &1.2421 & 1.1278 &0.5325 & 0.0003 & 4.0604 & 4.1747 &-&   0.5484\\
% 125   &1.0225 & 0.7778 &0.7024 & 0.0003 & 4.2800 & 4.5247 &-&   0.8660\\
% \bottomrule
%\end{tabular}%}
%\label{tab:times}
%\end{table}

\subsubsection{The effect of wind power control through output caps}
\begin{figure}
\includegraphics[width=0.8\columnwidth]{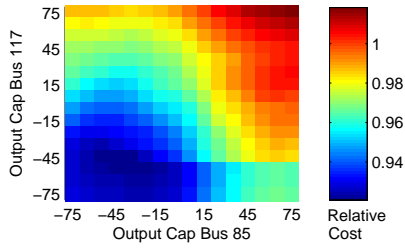}
\centering
\caption{Total cost of the case with caps on the output, relative to the case without any output caps. The output caps on the two different buses are defined relative to the mean output, and are varied from -75 MW to +75 MW. }
%The lowest cost is achieved when both wind power plants have an output cap 45 MW below the mean. The white region indicates infeasibility.}
\label{fig_cost_with_cap}
\end{figure}

In this section, we investigate the benefits of output cap control on wind power plants. We enable output cap control for the wind power plants at bus $85$ and bus $117$, since these have particularly large values of means and standard deviation, and allow the rest of the wind power plants to provide reserves according to Eq. \eqref{deltap_output}. We compare this with the case when no wind power plants have output caps, but all are able to provide reserves (including bus $85$ and $117$).
Figure~\ref{fig_cost_with_cap} shows the cost with cap control relative to the case without cap control, for various values of the cap threshold $\bar{\omega}_i$ in Eq.~(\ref{cap_output}).
Around the cap threshold values of $(-45 MW,-45 MW)$, the total cost is minimized, and is significantly lower ($\sim -6 \%$) than the case without caps but with reserves. The cost increases as we move away from this band.

\begin{figure}
\includegraphics[width=0.9\columnwidth]{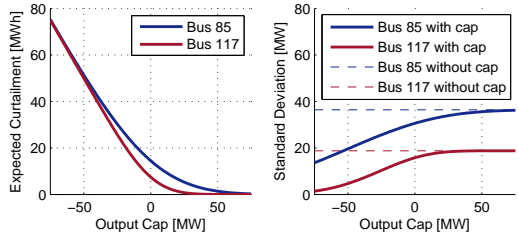}
\centering
\caption{The upper and lower plots shows the expected curtailed wind energy and and the standard deviation of the wind power output at bus 117 and bus 85, respectively. }
%While the expected curtailed energy decreases as the cap increases, the variability of the wind power output is higher at higher output caps.}
\label{fig_ExpCur_StdDev}
\end{figure}

To explain the trend in cost, we first investigate the effect of the cap threshold on the power output of the two wind power plants on bus $85$ and $117$.
In Figure~\ref{fig_ExpCur_StdDev}, the expected amount of wind power curtailment (left) and the standard deviation (right) for each power plant is shown as a function of the output cap on the left.
%In Figure~\ref{fig_ExpCur_StdDev}, the expected amount of wind power curtailment for each power plant is shown as a function of the output cap on the left. The behavior of the standard deviation of the wind power output as a function of the output cap is shown on the right. 
Note that, while the output caps are the same for both wind power plants, their nominal standard deviations (without caps) are different, which leads to different values for the expected curtailment and the standard deviation.
When the cap thresholds are lowered, the amount of utilized wind power drops due to higher expected curtailment. However, the standard deviation of the wind power output also drops significantly, which reduces wind power variability and thus the requirement for reserves.

%With lower output caps, the expected wind power curtailment increases. However, as opposed to $\Delta P$ control where 1 MW reduction of the mean $v$ leads to 1 MW curtailment, the curtailment effect of the output cap is non-linear.
%While the expected curtailment quantifies the lost energy, the standard deviation is a measure of the wind power output variability. For very high values of $\bar{\omega}$, the effect of capping is minimal, and the standard deviation is essentially the same as that for the uncapped power output. As we decrease $\bar{\omega}$ the standard deviation decreases towards zero. For low values of $\bar{\omega}$, the output of the wind power plant is close to constant.

%When the cap thresholds are lowered, the amount of utilized wind power drops due to higher expected curtailment. However, the uncertainty in wind power also drops significantly, which reduces requirement for reserves. 
%This phenomenon is quite clear for low values of the threshold.

\begin{figure}
\includegraphics[width=1.0\columnwidth]{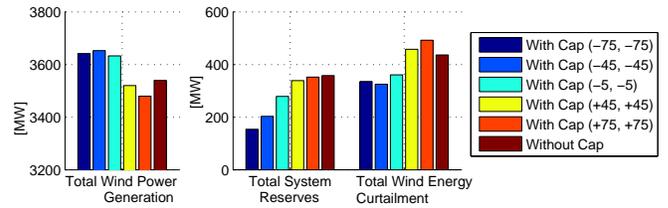}
\centering
\caption{Comparison of the cases with and without caps on the wind power output at bus 117 and bus 85.
%In the case with caps, both wind power plants are capped at -45 MW.
From left to right, the total generation by all wind power plants, the total amount of reserves in the system and the expected total wind energy curtailment are shown.}
%In the case with wind power caps, we have a lower requirement for reserves due to lower variability of the wind power production. This means that less wind power energy must be curtailed to provide reserves, which increases the total wind power generation.}
\label{fig_comparison}
\end{figure}

To investigate how this trade-off influences the cost in more detail, we look at the total wind power used for generation, the total system reserves and the total wind energy curtailment (which includes both ``wasted" energy and energy curtailed to provide reserves). These values are plotted for different pairs of output caps in Figure~\ref{fig_comparison}.

With low output caps $(-75,-75)$, higher wind curtailment outweighs the benefits of reduced reserve requirements compared to case with the optimal output caps $(-45,-45)$.
As the output caps increase, the total system reserve requirements also increase, and reaches the highest value for the case without output caps.
At the same time, we observe a decrease in the use of wind power for nominal energy production, since more of it must be allocated for reserves.

Th lowest utilization of wind for energy production happens with output caps of $(+75,+75)$, which is observed to be the point with highest cost in Figure~\ref{fig_cost_with_cap}. In this case, high wind power variability and no possibility of procuring wind power reserves at buses $85$ and $117$ leads to an increased curtailment of the mean wind power $v$ at those buses.

In the case with output caps $(+45,+45)$, the total wind power generation is lower than in the case without output caps. Still, the total cost of the solution is approximately 1\% lower with caps $(+45,+45)$ than without caps as observed in Figure~\ref{fig_cost_with_cap}. This is because the output caps reduces variability and thus the risk of overload on some critical transmission lines, which allows for better utilization of the transmission capacity and dispatch of cheaper conventional generators.

%\begin{figure}
%\includegraphics[width=0.9\columnwidth]{Fig/Comparison_WithWithout.eps}
%\centering
%\caption{Comparison of the cases with and without caps on the wind power output at bus 117 and bus 85. The upper plot shows the total generation by all wind power plants, the total amount of reserves in the system and the total wind power curtailment. The lower plots shows the expected curtailed energy and and the standard deviation of the wind power output at bus 117 and bus 85, respectively. We observe that the variability of the wind power plants are much lower in the case with output caps, which leads to a lower requirement for reserves. Although the total wind power curtailment is larger at bus 117 and bus 85 in the case with caps, the decreased requirement for reserves lead to a lower curtailment of wind power energy. Thus, the total produced energy from the wind power increases, reducing the cost.}
%\label{fig_comparison}
%\end{figure}
\section{Summary and conclusions}
\label{sec:Conclusion}
In this paper, we investigate two types of wind power control capabilities, $\Delta P$ control which reduces the mean power output by a constant value and allows the wind power plants to provide reserves, and output cap control which enforces a hard threshold on the total power output. 
%We analyze the benefits of using the curtailed wind power for reserve provision. 
The corresponding control and reserve models are developed and incorporated into an Optimal Power Flow formulation with Weighted Chance Constrains, which allows for controlling the risk of overloads due to wind fluctuations. Based on a case study on the IEEE 118 bus system, we observe that using $\Delta P$ control with wind reserves provides substantial cost benefits. For wind power plants with high variability, the output cap control is shown to outperform  $\Delta P$ control with reserves.

% The 'triggered' command can be changed if desired:



%\bibliographystyle{IEEEtran}
%\bibliography{yourbibliofile.bib(absolute or relative path)}{}


\bibliographystyle{IEEEtran}
\bibliography{20151001_bib_LR_v2,20151004_bib_SM}

\end{document}